\documentclass{amsart}%
\usepackage{amsfonts}
\usepackage{amsmath}
\usepackage{amssymb}
\usepackage{enumerate}
\usepackage{graphicx}%
\newtheorem{theorem}{Theorem}
\theoremstyle{plain}
\newtheorem{acknowledgement}{Acknowledgement}

\newtheorem{remark}{Remark}

\numberwithin{equation}{section}
\begin{document}
\title[Positivity of Ricci curvature]{Positivity of Ricci curvature under the K\"{a}hler--Ricci flow}
\author{Dan Knopf}
\address{University of Texas at Austin}
\email{danknopf@math.utexas.edu}
\urladdr{http://www.ma.utexas.edu/\symbol{126}danknopf/}
\thanks{Partially supported by NSF grant DMS-0328233.}
\date{\today}
\maketitle

\section{Introduction}

In analyzing a geometric evolution equation, it is of critical importance to
determine which curvature conditions are preserved or attained under the flow.
This is especially true for the Ricci and K\"{a}hler--Ricci flows, because of
the fundamental ways in which these flows can exploit curvature conditions to
reveal topological and canonical geometric properties of their underlying manifolds.

\smallskip

An \emph{invariant cone }in the space of curvature operators is one that is
preserved by a flow. For Ricci flow, the condition $R\geq0$ is preserved in
all dimensions, while the condition $R\leq0$ is preserved only in real
dimension two. Positive curvature operator is preserved in all dimensions
\cite{Ham86}, but positive sectional curvature is not preserved in dimensions
four and above. The known counterexamples, constructed recently by Ni
\cite{Ni}, are non-product metrics on the total space of the tangent bundles
$\mathbb{R}^{n}\hookrightarrow T\mathcal{S}^{n}\twoheadrightarrow
\mathcal{S}^{n}$ for $n\geq2$. However, positive sectional curvature and
positive Ricci curvature are preserved in dimension three \cite{Ham82}. In
higher dimensions, Huisken \cite{Huisken85} and Nishikawa \cite{Nish1, Nish2}
have demonstrated the invariance of certain sets defined by curvature pinching
conditions. Margerin \cite{Margerin} improved Huisken's constant, thereby
establishing a sharp invariant set in dimension four. Also in dimension four,
H.~Chen proved that $2$-positivity of the curvature operator is preserved
\cite{Chen}, while Hamilton later proved that positive isotropic curvature
constitutes an invariant cone \cite{Ham97}.

An \emph{attractive cone }in the space of curvature operators is one that is
entered asymptotically at a finite-time singularity. For example, the estimate
$R\geq R_{\min}(0)$ implies that the half-space $R\geq0$ constitutes an
attractive cone in all dimensions. Attractive cones place valuable
restrictions on which singularity models may appear. Important examples are
the pinching estimates for $3$-manifolds proved independently by Ivey
\cite{Ivey} and Hamilton \cite{Ham95}. These estimates imply that any rescaled
limit formed at a finite-time singularity in dimension three has nonnegative
sectional curvature. In turn, this fact is of fundamental importance in
Perelman's use of $\kappa$-solutions \cite{P1, P2} in his recent spectacular
progress toward the Geometrization Conjecture.

For K\"{a}hler--Ricci flow, Bando \cite{Bando} and Mok \cite{Mok} proved that
positive bisectional curvature defines an invariant cone. Positive Ricci
curvature is a more problematic condition (except of course in complex
dimension one). It is not expected to be an attractive cone. Indeed, certain
K\"{a}hler--Ricci solitons constructed by Feldman, Ilmanen, and the author
\cite{FIK} have Ricci curvature of mixed sign, yet are strongly conjectured to
occur as rescaled limits of singularities encountered by perturbed metrics on
the Calabi manifolds $\mathbb{CP}^{1}\hookrightarrow F^{k}\twoheadrightarrow
\mathbb{CP}^{n-1}$. (See Example 2.2 of \cite{FIK}.)

On the other hand, it would be very useful indeed if positive Ricci curvature
constituted an invariant cone for K\"{a}hler--Ricci flow. For example, parts
of the proofs of some recent important results of X.-X.~Chen and Tian for
K\"{a}hler--Ricci flow \cite{CT02, CT04} only require positive Ricci
curvature, but the authors impose the much stronger assumption of positive
bisectional curvature to be sure that positive Ricci curvature will be
preserved. (These papers also assume the existence of a K\"{a}hler--Einstein
metric on the underlying manifold. See in particular Remarks 1.5 -- 1.7 and
Question 9.3 of \cite{CT04}.) While positive Ricci curvature is not expected
to be preserved in high complex dimensions, there were until now promising
reasons to hope that it might be preserved in complex dimension two.
Unpublished work of Cao and Hamilton proves that positive orthogonal
bisectional curvature $B(x,x^{\perp})\geq0$ is preserved under
K\"{a}hler--Ricci flow and that its presence ensures that positive Ricci
curvature is also preserved \cite{CH}. Very recent results of Phong and Sturm
prove that positive Ricci curvature is preserved on a compact complex surface
with the property that the sum of any two eigenvalues of its traceless
curvature operator on traceless $(1,1)$-forms is non-negative \cite{PS1}. (See
also \cite{PS2}.)

\smallskip

The main purpose of this short note is to show that positive Ricci curvature
is \emph{not }preserved for arbitrary complex surfaces. We construct complete
K\"{a}hler metrics of bounded curvature and non-negative Ricci curvature whose
K\"{a}hler--Ricci evolutions immediately acquire Ricci curvature of mixed
sign. Let $L_{-1}^{2}$ denote $\mathbb{C}^{2}$ blown up at the origin. (See
below for an explanation of the notation.) In this paper, we will prove the following:

\begin{theorem}
\label{Main}There exists a one-parameter family (modulo scaling) of K\"{a}hler
metrics on $L_{-1}^{2}$ that are complete, of bounded curvature, and of
non-negative Ricci curvature, but which immediately develop mixed Ricci
curvature when evolved by the K\"{a}hler--Ricci flow.
\end{theorem}

In fact, Theorem \ref{Main} is a special case of a more general result. For
$n\geq2$ and $k\in\mathbb{N}$, let $L_{-k}^{n}$ denote the total space of the
holomorphic line bundle $\mathbb{C}\hookrightarrow L_{-k}^{n}%
\twoheadrightarrow\mathbb{CP}^{n-1}$ characterized by the equation
$\left\langle c_{1},[\Sigma]\right\rangle =-k$, where $c_{1}$ is the first
Chern class of the bundle and $\Sigma\approx\mathbb{CP}^{1}$ is a positively
oriented generator of $H_{2}(\mathbb{CP}^{n-1};\mathbb{Z})$. Since the
complement of the zero-section of $L_{-k}^{n}$ is biholomorphic to
$(\mathbb{C}^{n}\backslash\{0\})/\mathbb{Z}_{k}$, one may construct
$L_{-k}^{n}$ by gluing a $\mathbb{CP}^{n-1}$ onto $(\mathbb{C}^{n}%
\backslash\{0\})/\mathbb{Z}_{k}$ at the origin. Notice that $L_{-1}^{n}$ is
simply the tautological line bundle, i.e.~$\mathbb{C}^{n}$ blown up at the origin.

\begin{theorem}
\label{Full}For each complex dimension $n\geq2$ and all integers $k\in
\lbrack1,n-1]$, there exists a one-parameter family (modulo scaling)  of
K\"{a}hler metrics on $L_{-k}^{n}$ that are complete, of bounded curvature,
and of non-negative Ricci curvature, but which immediately develop mixed Ricci
curvature when evolved by the K\"{a}hler--Ricci flow.
\end{theorem}

The remainder of this paper is organized as follows. Since the critical $n=2$
case is most interesting, we begin with a detailed proof of Theorem
\ref{Main}. In Section \ref{Review} we provide, for the convenience of the
reader, a review of $\operatorname*{U}(2)$-invariant metrics on $\mathbb{C}%
^{2}\backslash\{0\}$. In Section \ref{Construction}, we construct initial
metrics on $L_{-1}^{2}$. Then in Section \ref{Evolve}, we finish the proof of
Theorem \ref{Main}. Finally, in Section \ref{General}, we indicate the
straightforward generalizations by which one proves Theorem \ref{Full} for all
complex dimensions $n\geq2$.

\begin{acknowledgement}
The author warmly thanks Huai-Dong Cao for encouraging him to work on this problem.
\end{acknowledgement}

\begin{remark}
Like Ni's counterexamples for non-negative sectional curvature \cite{Ni}, the
manifolds constructed here are noncompact. Thus it is still an open question
of significant interest to find compact counterexamples for which non-negative
Ricci (or sectional) curvatures are not preserved under K\"{a}hler--Ricci flow.
\end{remark}

\section{$\operatorname*{U}(2)$-invariant K\"{a}hler metrics\label{Review}}

In this section, we review the construction of $\operatorname*{U}%
(2)$-invariant K\"{a}hler metrics on $\mathbb{C}^{2}\backslash\{0\}$. Such
metrics also appear in \cite{Cao96, Cao97} and \cite{FIK}. Consider
$\mathbb{C}^{2}\backslash\{0\}$ with complex coordinates $(z_{1},z_{2})$.
Define%
\[
u=|z_{1}|^{2},\qquad v=|z_{2}|^{2},\qquad\text{and}\qquad w=u+v.
\]
Any function $P:\mathbb{C}^{2}\backslash\{0\}\rightarrow\mathbb{R}$ depending
only on $w$ may be written as $P(r)$, where%
\[
r=\log w.
\]
Anticipating that we will later abuse notation by regarding $P$ as a function
of time as well, we shall use subscripts to denote (partial) derivatives with
respect to $r$. Define
\begin{equation}
\varphi(r)=P_{r}(r).
\end{equation}
Then $P$ is the K\"{a}hler potential for a metric $g$ if and only if the
derivatives $\varphi$ and $\varphi_{r}$ are everywhere positive. In this case,
it is straightforward to compute that%
\begin{equation}
g_{i\bar{\jmath}}=e^{-r}\varphi\delta_{ij}+e^{-2r}(\varphi_{r}-\varphi)\bar
{z}_{i}z_{j}\label{g}%
\end{equation}
and%
\begin{equation}
g^{i\bar{\jmath}}=e^{r}\varphi^{-1}\delta^{ij}+(\varphi_{r}^{-1}-\varphi
^{-1})z^{i}\bar{z}^{j}.\label{g-inverse}%
\end{equation}
With respect to the standard bases $(\partial/\partial z_{1},\partial/\partial
z_{2})$ and $(dz_{1},dz_{2})$, one can therefore regard $g$ and $g^{-1}$ as
the matrices%
\[%
\begin{pmatrix}
g_{1\bar{1}} & g_{1\bar{2}}\\
g_{2\bar{1}} & g_{2\bar{2}}%
\end{pmatrix}
=\frac{1}{w^{2}}%
\begin{pmatrix}
v\varphi+u\varphi_{r} & (\varphi_{r}-\varphi)\bar{z}_{1}z_{2}\\
(\varphi_{r}-\varphi)z_{1}\bar{z}_{2} & u\varphi+v\varphi_{r}%
\end{pmatrix}
\]
and%
\[%
\begin{pmatrix}
g^{1\bar{1}} & g^{2\bar{1}}\\
g^{1\bar{2}} & g^{2\bar{2}}%
\end{pmatrix}
=\frac{1}{\varphi\varphi_{r}}%
\begin{pmatrix}
u\varphi+v\varphi_{r} & (\varphi-\varphi_{r})\bar{z}_{1}z_{2}\\
(\varphi-\varphi_{r})z_{1}\bar{z}_{2} & v\varphi+u\varphi_{r}%
\end{pmatrix}
\]
respectively. In the standard coordinate chart, define%
\begin{equation}
G=-\log\det g=2r-\log\varphi-\log\varphi_{r}\label{G}%
\end{equation}
and%
\begin{equation}
\psi=G_{r}=2-\frac{\varphi_{r}}{\varphi}-\frac{\varphi_{rr}}{\varphi_{r}%
}.\label{psi}%
\end{equation}
Then the complex Ricci tensor $\operatorname*{Rc}=R_{i\bar{\jmath}}%
\,dz_{i}\,d\bar{z}_{j}$ may be represented in the same basis as%
\begin{equation}%
\begin{pmatrix}
R_{1\bar{1}} & R_{1\bar{2}}\\
R_{2\bar{1}} & R_{2\bar{2}}%
\end{pmatrix}
=\frac{1}{w^{2}}%
\begin{pmatrix}
v\psi+u\psi_{r} & (\psi_{r}-\psi)\bar{z}_{1}z_{2}\\
(\psi_{r}-\psi)z_{1}\bar{z}_{2} & u\psi+v\psi_{r}%
\end{pmatrix}
.\label{RicciMatrix}%
\end{equation}
Abusing terminology slightly, we will say $V$ is an eigenvector of Ricci
corresponding to the eigenvalue $\lambda$ if $R_{i\bar{\jmath}}V^{i}=\lambda
g_{i\bar{\jmath}}V^{i}$. Thus understood, the eigenvalues of
$\operatorname*{Rc}$ are%
\[
\lambda_{1}=\frac{\psi}{w}\qquad\text{and}\qquad\lambda_{2}=\frac{\psi_{r}}%
{w}.
\]

The Levi-Civita connection of $g$ is determined by the Christoffel symbols%
\begin{align*}
\Gamma_{11}^{1}  &  =\frac{\bar{z}_{1}}{w^{2}}\left[  u\frac{\varphi_{rr}%
}{\varphi_{r}}+2v\frac{\varphi_{r}}{\varphi}+u-2w\right] \\
\Gamma_{12}^{1}  &  =\frac{\bar{z}_{2}}{w^{2}}\left[  u\frac{\varphi_{rr}%
}{\varphi_{r}}+(v-u)\frac{\varphi_{r}}{\varphi}-v\right] \\
\Gamma_{22}^{1}  &  =\frac{z_{1}\bar{z}_{2}^{2}}{w^{2}}\left[  \frac
{\varphi_{rr}}{\varphi_{r}}-2\frac{\varphi_{r}}{\varphi}+1\right] \\
\Gamma_{11}^{2}  &  =\frac{\bar{z}_{1}^{2}z_{2}}{w^{2}}\left[  \frac
{\varphi_{rr}}{\varphi_{r}}-2\frac{\varphi_{r}}{\varphi}+1\right] \\
\Gamma_{12}^{2}  &  =\frac{\bar{z}_{1}}{w^{2}}\left[  v\frac{\varphi_{rr}%
}{\varphi_{r}}+(u-v)\frac{\varphi_{r}}{\varphi}-u\right] \\
\Gamma_{22}^{2}  &  =\frac{\bar{z}_{2}}{w^{2}}\left[  v\frac{\varphi_{rr}%
}{\varphi_{r}}+2u\frac{\varphi_{r}}{\varphi}+v-2w\right]
\end{align*}
and their conjugates. The complex Riemann tensor is determined by%
\[
R_{i\bar{\jmath}k\bar{\ell}}=-g_{m\bar{\ell}}\bar{\partial}_{j}\Gamma_{ik}%
^{m}=-\partial_{i}\bar{\partial}_{j}g_{k\bar{\ell}}+g^{p\bar{q}}\partial
_{i}g_{k\bar{q}}\bar{\partial}_{j}g_{p\bar{\ell}}.
\]
A straightforward computation shows that%
\begin{align*}
R_{i\bar{\jmath}k\bar{\ell}}  &  =e^{-4r}\left[  -\varphi_{rrr}+4\varphi
_{rr}-2\varphi_{r}+2\varphi-4\frac{\varphi_{r}^{2}}{\varphi}+\frac
{\varphi_{rr}^{2}}{\varphi_{r}}\right]  \bar{z}_{i}z_{j}\bar{z}_{k}z_{\ell}\\
&  +e^{-3r}\left[  -\varphi_{rr}+\varphi_{r}-\varphi+\frac{\varphi_{r}^{2}%
}{\varphi}\right]  (\bar{z}_{i}z_{j}\delta_{k\ell}+\bar{z}_{i}\delta
_{jk}z_{\ell}+\delta_{ij}\bar{z}_{k}z_{\ell}+\delta_{i\ell}z_{j}\bar{z}_{k})\\
&  +e^{-2r}\left[  -\varphi_{r}+\varphi\right]  (\delta_{ij}\delta_{k\ell
}+\delta_{i\ell}\delta_{jk}).
\end{align*}
To put this into a more useful form, we follow \cite{Cao96}. Since the metric
is $\operatorname*{U}(2)$-invariant, there is no loss of generality in
evaluating the curvature at $(\zeta,0)$, where $\zeta\neq0$ is arbitrary.
Furthermore, one may assume that the complex tangent space at $(\zeta,0)$ is
spanned by an orthonormal basis of complex vectors $x,y$ such that $x^{2}=0$.
Then one can readily compute the bisectional curvatures $B(x,x)=R(x,\bar
{x},x,\bar{x})$, $B(x,y)=R(x,\bar{x},y,\bar{y})$, and $B(y,y)=R(y,\bar
{y},y,\bar{y})$ at $(\zeta,0)$, obtaining%
\begin{align*}
B(x,x)|_{(\zeta,0)}  &  =\frac{1}{|\zeta|^{4}}(\frac{\varphi_{rr}^{2}}%
{\varphi_{r}}-\varphi_{rrr})|x^{1}|^{4}\\
B(x,y)|_{(\zeta,0)}  &  =\frac{1}{|\zeta|^{4}}\left[  (\frac{\varphi_{rr}^{2}%
}{\varphi_{r}}-\varphi_{rrr})|x^{1}|^{2}|y^{1}|^{2}+(\frac{\varphi_{r}^{2}%
}{\varphi}-\varphi_{rr})|x^{1}|^{2}|y^{2}|^{2}\right] \\
B(y,y)|_{(\zeta,0)}  &  =\frac{1}{|\zeta|^{4}}\left[  (\frac{\varphi_{rr}^{2}%
}{\varphi_{r}}-\varphi_{rrr})|y^{1}|^{4}+4(\frac{\varphi_{r}^{2}}{\varphi
}-\varphi_{rr})|y^{1}|^{2}|y^{2}|^{2}+2(\varphi-\varphi_{r})|y^{2}%
|^{4}\right]  .
\end{align*}
Observe in particular that the choices%
\begin{equation}
X|_{(\zeta,0)}=\frac{|\zeta|}{\sqrt{\varphi_{r}}}\,\frac{\partial}{\partial
z_{1}}\qquad\text{and}\qquad Y|_{(\zeta,0)}=\frac{|\zeta|}{\sqrt{\varphi}%
}\,\frac{\partial}{\partial z_{2}} \label{XY}%
\end{equation}
allow one to recover (\ref{RicciMatrix}) in the form%
\begin{align*}
\operatorname*{Rc}(X,\bar{X})|_{(\zeta,0)}  &  =B(X,X)+B(X,Y)=\frac{\psi_{r}%
}{\varphi_{r}}\\
\operatorname*{Rc}(Y,\bar{Y})|_{(\zeta,0)}  &  =B(Y,Y)+B(X,Y)=\frac{\psi
}{\varphi}.
\end{align*}

\section{An initial metric of positive Ricci curvature\label{Construction}}

We now want to investigate whether there exists a $\operatorname*{U}%
(2)$-invariant K\"{a}hler potential $P$ on $\mathbb{C}^{2}\backslash\{0\}$
with the following global properties:

\begin{enumerate}
\item \label{g1}$\varphi>0$ everywhere;

\item \label{g2}$\varphi_{r}>0$ everywhere;

\item \label{R1}$\psi>0$ everywhere;

\item \label{R2}$\psi_{r}\geq0$ everywhere;

\item \label{B}$g$ has bounded curvature;

\item \label{C1}$g$ extends smoothly to a complete metric as $r\rightarrow
+\infty$; and

\item \label{C2}$g$ extends smoothly to a complete metric as $r\rightarrow
-\infty$.
\end{enumerate}

Furthermore, we require that there exists at least one point $(\zeta
,0)\in\mathbb{C}^{2}\backslash\{0\}$ such that with respect to an orthonormal
basis $x,y$ for the complex tangent space at $(\zeta,0)$, the following local
properties hold:

\begin{enumerate}
[I]

\item $B(x,y)<0$; and

\item $B(x,x)=-B(x,y)$.
\end{enumerate}

It may seem foolhardy to expect a single function $P$ to satisfy all nine of
these requirements. However, let $a>0$ be a constant to be determined and
consider the second-order \textsc{ode}%
\begin{equation}
\left[  \log(\varphi\varphi_{r})\right]  _{r}=a.
\end{equation}
Its general solution is $\sqrt{2e^{ar+b}/a+c}$, where $b$ and $c$ are
arbitrary. Without loss of generality, one may eliminate $b$ by scaling,
yielding%
\begin{equation}
\varphi(r)=\sqrt{\frac{2}{a}e^{ar}+c}. \label{phi}%
\end{equation}
One then has
\begin{equation}
\varphi\varphi_{r}=e^{ar}, \label{phi-phi-r}%
\end{equation}
so that global properties (\ref{g1}) and (\ref{g2}) are satisfied whenever
$c\geq0$. Moreover,
\[
\psi=2-a
\]
everywhere. So properties (\ref{R1}) and (\ref{R2}) are satisfied whenever
$a\in(0,2)$.

To evaluate the bisectional curvatures of the resulting metric, it suffices to
consider the orthonormal basis (\ref{XY}) for the complex tangent space at
$(\zeta,0)$, where $\zeta\neq0$ is arbitrary. One finds that%
\begin{equation}
B(X,X)=\frac{ac}{\varphi^{3}}=-B(X,Y), \label{BXY}%
\end{equation}
so that local properties I and II will in fact be satisfied globally if $c>0$.
The choices $a\in(0,2)$ and $c>0$ also ensure that%
\[
B(Y,Y)=\frac{2}{\varphi^{3}}\left[  \left(  \frac{2}{a}-1\right)
e^{ar}+c\right]  >B(X,X).
\]
Because $\varphi(r)\rightarrow\sqrt{c}>0$ as $r\rightarrow-\infty$ and
$\varphi(r)=\mathcal{O}(e^{(a/2)r})$ as $r\rightarrow\infty$, the curvatures
are globally bounded, satisfying property (\ref{B}). To investigate
completeness at spatial infinity, observe that (\ref{phi-phi-r}) implies that
$\varphi_{r}=\mathcal{O}(e^{ar})$ as $r\rightarrow-\infty$ and that
$\varphi_{r}=\mathcal{O}(e^{(a/2)r})$ as $r\rightarrow\infty$. Since radial
paths emanating from the origin are geodesics and the length of such a path is
$\int_{-\infty}^{\infty}\sqrt{\varphi_{r}}\,dr$, it is easy to see that the
metric is complete as $r\rightarrow+\infty$, satisfying global property
(\ref{C1}).

On the other hand, the metric distance from an arbitrary point to $r=-\infty$
is finite. Our choice of $c>0$ makes it impossible to complete the metric
smoothly by gluing in a point at the origin. But Calabi's theorem
\cite{Calabi} tells us that $g$ will extend to a smooth K\"{a}hler metric on
$L_{-k}^{2}$ if $\varphi$ has an asymptotic expansion near $w=0$ of the form%
\[
\varphi=a_{0}+a_{1}w^{k}+a_{2}w^{2k}+\mathcal{O}(w^{3k}),
\]
with $a_{0}$ and $a_{1}$both positive. (See Section 4.1 of \cite{FIK} for an
exposition.) The function (\ref{phi}) will have such an expansion if $a$ is a
positive integer. Fortunately for us, there is a (unique) positive integer in
the interval $(0,2)$. Taking $a=1$ thus lets us satisfy global property
(\ref{C2}). Notice that for each $c>0$, the solution $\varphi$ then admits the
expansion
\[
\varphi=\allowbreak\sqrt{c}+\frac{1}{\sqrt{c}}w-\frac{1}{2c^{3/2}}w^{2}%
+\cdots.
\]
In particular, $2\pi\sqrt{c}$ is the area of the $\mathbb{CP}^{1}$ that
constitutes the zero-section of the bundle.

\begin{remark}
The K\"{a}hler manifolds $(L_{-1}^{2},g)$ we have constructed are asymptotic
at infinity to a K\"{a}hler cone $\Gamma$ on $\mathbb{C}^{2}\backslash\{0\}$
whose cone angle around the origin in each complex line is $\pi$ and whose
density ratio is $\operatorname*{Vol}_{\Gamma}(B_{1})/\operatorname*{Vol}%
_{\mathbb{C}^{2}}(B_{1})=1/4$.
\end{remark}

\section{Effect of the K\"{a}hler--Ricci flow\label{Evolve}}

We now consider the K\"{a}hler--Ricci evolution $(L_{-1}^{2},g(t))$, where we
take $P(r,0)$ to be the K\"{a}hler potential $P(r)$ constructed above. Again
choose an arbitrary point $(\zeta,0)\in\mathbb{C}^{2}\backslash\{0\}$. Because%
\[
g|_{(\zeta,0)}=\frac{1}{|\zeta|^{2}}%
\begin{pmatrix}
\varphi_{r} & 0\\
0 & \varphi
\end{pmatrix}
\]
and%
\[
\operatorname*{Rc}|_{(\zeta,0)}=\frac{1}{|\zeta|^{2}}%
\begin{pmatrix}
\psi_{r} & 0\\
0 & \psi
\end{pmatrix}
\]
in the standard basis, the K\"{a}hler--Ricci flow%
\[
\frac{\partial}{\partial t}g=-\operatorname*{Rc}%
\]
will be satisfied if and only if $\varphi$ evolves by the \textsc{pde}
$\varphi_{t}=-\psi$; to wit,%
\begin{equation}
\varphi_{t}=\frac{\varphi_{rr}}{\varphi_{r}}+\frac{\varphi_{r}}{\varphi
}-2.\label{phi-t}%
\end{equation}
Since $\partial/\partial r$ and $\partial/\partial t$ commute, it follows that
the Ricci potential $\psi$ must evolve by the \textsc{pde}%
\begin{equation}
\psi_{t}=\frac{\psi_{rr}}{\varphi_{r}}-\frac{\varphi_{rr}\psi_{r}}{\varphi
_{r}^{2}}+\frac{\psi_{r}}{\varphi}-\frac{\varphi_{r}\psi}{\varphi^{2}%
}.\label{psi-t}%
\end{equation}
Since $\psi(\cdot,0)\equiv1$, this reduces to%
\[
\left.  \frac{\partial}{\partial t}\psi\right\vert _{t=0}=-\frac{\varphi_{r}%
}{\varphi^{2}}%
\]
at the initial time. A further computation then shows that%
\begin{equation}
\left.  \frac{\partial}{\partial t}\psi_{r}\right\vert _{t=0}=\frac{e^{r}%
}{\varphi^{5}}(e^{r}-c),\label{key}%
\end{equation}
which is strictly negative for all $r<\log c$. Since $\psi_{r}(\cdot
,0)\equiv0$, the complex Ricci tensor must acquire a negative eigenvalue
$\lambda_{2}<0$ at all such points for small times $t>0$. This completes the
proof of Theorem \ref{Main}.

\section{The general case: complex dimension $n\geq2$\label{General}}

If $P$ is a $\operatorname*{U}(n)$-invariant K\"{a}hler potential on
$\mathbb{C}^{n}\backslash\{0\}$, formulas (\ref{g}) for $g$ and
(\ref{g-inverse}) for $g^{-1}$ remain unchanged. Formula (\ref{G}) becomes%
\[
G=nr-\log\varphi-\log\varphi_{r},
\]
whence (\ref{psi}) is replaced by%
\[
\psi=G_{r}=n-\frac{\varphi_{r}}{\varphi}-\frac{\varphi_{rr}}{\varphi_{r}}.
\]
One calculates easily that%
\[
R_{i\bar{\jmath}}=e^{-r}\psi\delta_{ij}+e^{-2r}(\psi_{r}-\psi)\bar{z}_{i}%
z_{j}.
\]
To display the eigenvalues of $\operatorname*{Rc}$, one may without loss of
generality fix $Z=(\zeta,0,\ldots,0)\in\mathbb{C}^{n}\backslash\{0\}$ with
$\zeta\neq0$ arbitrary. Then%
\[
g|_{Z}=\frac{1}{|\zeta|^{2}}%
\begin{bmatrix}
\varphi_{r} &  &  & \\
& \varphi &  & \\
&  & \ddots & \\
&  &  & \varphi
\end{bmatrix}
\]
and%
\[
\operatorname*{Rc}|_{Z}=\frac{1}{|\zeta|^{2}}%
\begin{bmatrix}
\psi_{r} &  &  & \\
& \psi &  & \\
&  & \ddots & \\
&  &  & \psi
\end{bmatrix}
.
\]

If $\varphi$ is defined by (\ref{phi}) with $a>0$, one finds that%
\[
\psi\equiv n-a,
\]
hence that $\psi_{r}\equiv0$. Equation (\ref{phi-t}) becomes%
\[
\varphi_{t}=\frac{\varphi_{rr}}{\varphi_{r}}+\frac{\varphi_{r}}{\varphi}-n,
\]
so that (\ref{psi-t}) remains unchanged. The analog of (\ref{key}) is%
\[
\left.  \frac{\partial}{\partial t}\psi_{r}\right\vert _{t=0}=(n-a)\frac
{e^{ar}}{\varphi^{5}}(e^{ar}-ac),
\]
which is strictly negative for all $a\in(0,n)$ and $c>0$ whenever $r<\frac
{1}{a}\log(ac)$. If we take $a\in(0,n)$ to be an integer, then $\varphi$ will
have an asymptotic expansion near the origin of the form%
\[
\varphi=\sqrt{c}+\frac{1}{a\sqrt{c}}w^{a}-\frac{1}{2a^{2}c^{3/2}}w^{2a}%
+\cdots.
\]
Hence by Calabi's theorem, $g$ will extend to a smooth metric on $L_{-a}^{n}$
for all such $a$. It is easy to check that $g$ is complete and asymptotic to a
K\"{a}hler cone at infinity, hence of bounded curvature. This completes the
proof of Theorem \ref{Full}. 

\begin{remark}
For each $1\leq k\leq n-1$, the K\"{a}hler manifolds $(L_{-k}^{n},g)$ we have
constructed are asymptotic at infinity to a K\"{a}hler cone $\Gamma_{k}$ on
$\mathbb{C}^{n}\backslash\{0\}$ whose cone angle around the origin is $\pi$
and whose density ratio is $\operatorname*{Vol}_{\Gamma_{k}}(B_{1}%
)/\operatorname*{Vol}_{\mathbb{C}^{n}}(B_{1})=k^{1-n}/2^{n}$.
\end{remark}

\end{document}